\documentclass[12pt]{amsart}

\usepackage{amssymb,mathrsfs,booktabs,url}
\usepackage{xspace,verbatim}
\usepackage[noadjust]{cite}
\usepackage[line,color,poly,all]{xy}
\usepackage[breaklinks=true]{hyperref}

\renewcommand{\O}{\mathcal{O}}
\newcommand{\Gal}{\operatorname{Gal}}

\newcommand{\p}{\mathfrak{p}}

\newcommand{\Z}{\mathbb{Z}}
\newcommand{\F}{\mathbb{F}}
\newcommand{\R}{\mathbb{R}}

\newcommand{\Q}{\mathbb{Q}}

\newcommand{\U}{\operatorname{U}}

\newcommand{\GL}{\operatorname{GL}}
\newcommand{\SL}{\operatorname{SL}}

\newcommand{\C}{\mathbb{C}}

\newcommand{\Res}{\operatorname{Res}}

\theoremstyle{definition}

\newcommand{\CC}{\mathbb{C}}
\newcommand{\QQ}{\mathbb{Q}}

\newcommand{\z}{\zeta}
\newcommand{\fn}{\mathfrak{n}}
\newcommand{\fp}{\mathfrak{p}}
\newcommand{\fq}{\mathfrak{q}}
\newcommand{\fr}{\mathfrak{r}}
\newcommand{\pphi}{\hat{\phi}} 
\newcommand{\vect}[1]{\begin{pmatrix} #1 \end{pmatrix}}%
\newcommand{\mat}[1]{\begin{pmatrix} #1 \end{pmatrix}}%
\newcommand{\Vor}{Vorono\"{\i}\xspace}
\DeclareMathOperator{\Norm}{Norm}
\DeclareMathOperator{\Tr}{Tr}

\newcommand{\bG}{\mathbf{G}}
\newcommand{\OO}{\O}
\newcommand{\fH}{\mathfrak{H}}
\newcommand{\fm}{\mathfrak{m}}
\newcommand{\cusp}{\operatorname{cusp}}
\newcommand{\V}{{\mathcal V}}

\newcommand{\cC}{\mathcal{C}}
\newcommand{\bcC}{\bar{\cC}}
\newcommand{\ZZ}{\Z}
\newcommand{\del}{\partial }

\newcommand{\vv}{\mathbf{v}}
\newcommand{\longisomto}{\overset{\sim}{\longrightarrow}}
\newcommand{\bs}{\backslash }
\DeclareMathOperator{\Eis}{Eis}
\DeclareMathOperator{\sgn}{sgn}
\DeclareMathOperator{\Size}{size}

\begin{document}

\title[Modular forms and elliptic curves]{Modular forms and elliptic curves over the field of fifth roots
of unity}
\author{Paul E. Gunnells}
\address{Department of Mathematics and Statistics\\University of Massachusetts\\Amherst, MA 01003-9305}
\email{gunnells@math.umass.edu}
\author{Farshid Hajir}
\address{Department of Mathematics and Statistics\\University of Massachusetts\\Amherst, MA 01003-9305}
\email{hajir@math.umass.edu}
\author{Dan Yasaki}
\address{Department of Mathematics and Statistics\\ 
University of North Carolina at Greensboro\\Greensboro, NC 27402-6170}
\email{d\_yasaki@uncg.edu}


\date{14 May 2010} 

\thanks{We thank Mark Watkins for writing the appendix to this paper.
PG thanks the NSF for support.  FH thanks the NSA for support.
This manuscript is submitted for publication with the understanding
that the United States government is authorized to produce and
distribute reprints. DY thanks UNC Greensboro for support through a
UNC Greensboro New Faculty Grant.}



\begin{abstract}
Let $F$ be the cyclotomic field of fifth roots of unity.  We
computationally investigate modularity of elliptic curves over $F$.
\end{abstract}

\maketitle

\section{Introduction}

Let $\zeta$ be a primitive fifth root of unity, and let
$F = \Q (\zeta)$.  In this paper we describe computational
work that investigates the modularity of elliptic curves over $F$.
Here by \emph{modularity} we mean that for a given elliptic curve $E$
over $F$ with conductor $\fn$ there should exist an automorphic form $f$
on $\GL_{2}$, also of conductor $\fn$, such that we have the equality of
partial $L$-functions $L^{S}(s, f) = L^{S}(s, E)$, where $S$ is a
finite set of places including those dividing $\fn$.  We are also
interested in checking a converse to this notion, which says that for
an appropriate automorphic form $f$ on $\GL_{2}$, there should exist
an elliptic curve $E/F$ again with matching of partial $L$-functions.
Our work is in the spirit of that of Cremona and his students
\cite{cremona2, crem.whitley,lingham, bygott} for complex quadratic
fields, and of Socrates--Whitehouse \cite{sw} and Demb{\'e}l{\'e}
\cite{dembele} for real quadratic fields.

Instead of working with automorphic forms, we work with the cohomology
of congruence subgroups of $\GL_{2} (\OO)$, where $\OO$ is the ring of
integers of $F$.  A main motivation for this is the Eichler--Shimura
isomorphism, which identifies the cohomology of subgroups of $\SL_{2}
(\Z)$ with spaces of modular forms.  More precisely, let $N\geq 1$ be 
an integer and let $\Gamma_{0} (N)\subset \SL_2(\Z )$ be the usual
congruence subgroup of matrices upper triangular mod $N$.  The group
cohomology $H^{*} (\Gamma_{0} (N); \C)$ is isomorphic to the
cohomology $H^{*} (\Gamma_{0} (N)\backslash \fH; \C)$, where $\fH$ is
the upper halfplane.  We have an isomorphism
\begin{equation}\label{eq:eichlershimura}
H^{1} (\Gamma_{0} (N); \C) \simeq S_{2}
(N)\oplus \overline{S}_{2} (N)\oplus \Eis_2 (N),
\end{equation}
where $S_{2} (N)$ is the space of weight two holomorphic cusp forms of
level $N$, the summand $\Eis_2(N)$ is the space of weight two
holomorphic Eisenstein series of level $N$, and the bar denotes complex
conjugation.

Moreover \eqref{eq:eichlershimura} is an isomorphism of Hecke modules:
there are Hecke operators defined on the cohomology $H^{1} (\Gamma_{0}
(N); \C)$ that parallel the usual operators defined on modular forms,
and the two actions respect the isomorphism.  This means that the
cohomology of $\Gamma_{0} (N)$ provides a concrete way to compute
with the modular forms of interest in the study of elliptic curves
over $\Q$.

Further motivation is provided by Franke's proof of Borel's conjecture
\cite{Fra}.  Franke's work shows that the cohomology of arithmetic
groups can always be computed in terms of certain automorphic forms.
Although the forms that occur in cohomology are a small subset of all
automorphic forms, they are widely believed to have deep connections
with arithmetic geometry.  In particular, let $\Gamma_{0} (\fn)\subset
\GL_2(\OO) $ be the congruence subgroup of matrices upper triangular
modulo $\fn$.  There is a subspace of the cohomology $H^{*}
(\Gamma_{0} (\fn); \C)$ called the \emph{cuspidal cohomology} that
corresponds to cuspidal automorphic forms.  This subspace provides a
natural place to realize the ``appropriate'' automorphic forms above.
Thus instead of defining what a ``weight $2$ modular form over $F$ of
level $\fn$'' means, we work with the cuspidal cohomology with trivial
coefficients of the congruence subgroup $\Gamma_{0} (\fn)$.

We now give an overview of the contents of this paper and summarize
our main results.  In \S\ref{s:gb} we give the geometric background of
our cohomology computations and describe the Hecke operators and how
they act on cohomology.  The next two sections give details about how
we performed the cohomology computations.  In \S\ref{s:rt} we explain
the explicit reduction theory we need for the group $\GL_{2} (\O)$,
and in \S\ref{s:ho} we discuss how we compute the action of the Hecke
operators on cohomology.  Next we turn to the elliptic curve side of
the story, and in \S\ref{s:ec} we examine various methods for writing
down elliptic curves over $F$.  Here the methods are more ad hoc than
on the cohomology side.  We describe the straightforward method of
searching ``in a box,'' and a trick using $S$-unit equations and the
Frey--Hellegouarch construction.  Finally in \S\ref{s:data} we present
our computational data.  We give tables of cohomology data, including
the levels where we found cuspidal cohomology and the dimensions, as
well as some eigenvalues of Hecke operators $T_{\fq}$ for a range of
primes $\fq$.  We then give ``motivic'' explanations for the cuspidal
cohomology classes with rational Hecke eigenvalues --- either by
identifying them as arising from weight $2$ modular forms on $\Q$ or
parallel weight $2$ Hilbert modular forms on $F^+=\Q(\sqrt{5})$, or by
finding elliptic curves over $F$ --- that apparently match the
eigenvalue data.
 
We were able to motivically account for every rational Hecke
eigenclass we computed.  All eigenclasses that appeared to come from
classes over $\Q$ and $F^{+}$ were found using tables computed by
Cremona \cite{cremona} and tables/software due to Demb\'el\'e
\cite{dembele}.  Of the eigenclasses that do not come from $\Q$ and
$F^{+}$, for all but one our searches found elliptic curves over
$\Q(\zeta_5)$ whose point counts matched the eigenvalue data.  We also
note that one rational eigenclass we found corresponds to a ``fake
elliptic curve'' in the sense of Cremona \cite{cremona.fake}.  Details
can be found in \S\ref{s:data}.  The only form we were unable to
account for occured at norm level $3641$.  After this
paper was first distributed, Mark Watkins conducted a successful
targeted search for the missing curve by modifying techniques of
Cremona--Lingham \cite{cl}.  We thank him for writing an appendix 
describing his result.

Conversely, within the range of our computations we were able to
cohomologically account for all the elliptic curves over $F$ that we
found.  That is, we found no elliptic curve over $F$ that was not
predicted by a rational Hecke eigenclass.

\subsection*{Acknowledgements} We thank Avner Ash, Kevin Buzzard, John
Cremona, and Lassina Demb\'el\'e for helpful conversations and
correspondence.  We especially thank Dinakar Ramakrishnan for
suggesting this project and for his encouragement.  Finally, we thank
Mark Watkins for finding the missing curve at norm level 3641 and for
writing the appendix.

\section{Geometric background}\label{s:gb}

\subsection{}\label{ss:caf} Let $\bG$ be the reductive $\Q$-group
$\Res _{F/\Q} (\GL_{2})$, where $\Res $ denotes restriction of
scalars.  We have $\bG(\Q) \simeq \GL_{2} (F)$.  Let $G=\bG(\R)$ be
the group of real points.  We have $G\simeq \GL_{2} (\C) \times
\GL_{2} (\C)$, where the two factors corresponding to the two
non-conjugate pairs of complex embeddings of $F$.  Let $K\simeq
\U(2)\times \U(2)$ be the maximal compact subgroup of $G$, and let
$A_{G}\simeq \C^{\times}$ be the identity component of the real
points of the maximal $\Q$-split torus in the center of $G$.  Fix an
ideal $\fn \subset \OO$, and let $\Gamma$ be the congruence subgroup
$\Gamma_{0} (\fn)$ defined in the introduction.

Let $X$ be the global symmetric space $G/A_{G}K$.  We have an
isomorphism 
\begin{equation}\label{eq:symspace}
X\simeq \fH_3\times \fH_3\times \R,
\end{equation}
where $\fH_3$ is hyperbolic $3$-space; thus $X$ is $7$-dimensional.

The space $X$ should be compared with the product of upper halfplanes
$\fH \times \fH$ one sees when studying Hilbert modular forms over
quadratic fields.  Indeed, if we were to work instead with $\bG ' =
R_{F/\Q} (\SL_2)$, the appropriate symmetric space would be
$\fH_{3}\times \fH_{3}$, which makes the analogy clear.  The extra
flat factor $\R$ in \eqref{eq:symspace} accounts for the difference
between the centers of $\GL_{2} (\OO)$ and $\SL_{2} (\OO)$.  As we
will see in \S\ref{s:rt}, it is much more convenient computationally
to work with $\GL_{2}$ instead of $\SL_{2}$.

\subsection{} We are interested in the complex group cohomology $H^{*}
(\Gamma ; \C)$, which can be identified with $H^{*} (\Gamma \backslash
X; \C)$.  As mentioned in the introduction, there is a precise way to
compute these cohomology spaces in terms of automorphic forms, and
there is a distinguished subspace $H^{*}_{\cusp}(\Gamma \backslash X;
\C)$ corresponding to the cuspidal automorphic forms.  We will not
make this explicit here, and instead refer to
\cite{li.schwermer.survey, harder.book, bw, borel.survey} for more
information.  Our goal now is to pin down exactly which cohomology
group we want to study.  In other words, which cohomology space $H^{i}
(\Gamma \backslash X; \C)$, where $0\leq i\leq 7$, plays the role of
$H^{1}$ of the modular curve?

First, although we a priori have cohomology in degrees $0$ to $7$, a
result of Borel--Serre \cite{BS} implies that $H^{7}$ vanishes
identically.  Moreover, standard computations from representation
theory (cf.~\cite{li.schwermer.survey}) show that $H^{i}_{\cusp}
(\Gamma \backslash X; \C)=0$ unless $2\leq i\leq 5$.  One also knows
that if a cuspform contributes to any of these degrees, it does to
all, and in essentially the same way.  For computational reasons it is
much easier to work with cohomology groups of higher degree, and so we
choose to work with $H^{5} (\Gamma \backslash X; \C)$.

\subsection{}
Next we consider the Hecke operators.  Let $\tilde{\Gamma}\subset \bG
(\Q)$ be the commensurator of $\Gamma$.  By definition
$\tilde{\Gamma}$ consists of all $g\in \bG (\Q)$ such that both
$\Gamma$ and $\Gamma^{g} := g^{-1}\Gamma g$ have finite index in
$\Gamma ' := \Gamma \cap \Gamma^g$.
The inclusions $\Gamma '\rightarrow \Gamma $
and $\Gamma '\rightarrow \Gamma^ g$ determine a diagram
\[
\vbox{\xymatrix{&{\Gamma '\backslash X}\ar[dr]^{t}\ar[dl]_{s}&\\
          {\Gamma \backslash X}&&{\Gamma \backslash X}}}
\]
Here $s(\Gamma 'x) = \Gamma x$ and $t$ is the composition of $\Gamma '
x \mapsto \Gamma ^{g}x$ with left multiplication by $g$.  This diagram
is the \emph{Hecke correspondence} associated to $g$.  It can be shown
that, up to isomorphism, the Hecke correspondence depends only on the
double coset $\Gamma g\Gamma $.

Because the maps $s$ and $t$ are proper, they induce a map on
cohomology:
\[
t_{*}s^{*}\colon H^{*}(\Gamma \backslash X;\Z
)\rightarrow H^{*}(\Gamma \backslash X;\Z). 
\]
We denote the induced map by $T_{g}$ and call it the \emph{Hecke
operator} associated to $g$. 

In our application we consider $g$ of the form $(\begin{smallmatrix}1&0\\
0&a
\end{smallmatrix}) $, where $a$ is a generator of any prime ideal $\fq$
coprime with $\fn$ (note every ideal in $\OO$ is principal since $F$
has class number $1$). Thus we are led to the main computational issue
on the modular side: for each $\fn$ compute the space $H^{5}_{\cusp}
(\Gamma_{0} (\fn)\backslash X; \C)$ together with the action of the
Hecke operators 
\[\bigl\{T_{\fq}\bigm| \text{$\fq $ prime, $\fq\nmid \fn$} \bigr\}.\]

\section{Reduction theory}\label{s:rt}

\subsection{} In this section we explain the connection between our
symmetric space $X$ and a cone of Hermitian forms.  This connection is
exactly the reason we prefer to work with $\GL_{2} (\OO)$ instead of
$\SL_{2} (\OO)$.  Let $\iota=(\iota_1,\iota_2)$ denote the (non
complex conjugate) embeddings
\[\iota : F \to \CC \times \CC\]
given by sending $\z$ to $(\z,\z^3)$.  We abbreviate the
second embedding by $\cdot'$, and for $\alpha\in F$ write $(\alpha,
\alpha')$ for $\iota(\alpha)$.

First let $V$ be the real vector space of $2\times 2$ Hermitian
matrices over $\C$. Let $C\subset V$ be the cone of positive-definite
Hermitian matrices.  The cone $C$ is preserved by homotheties (scaling
by $\R_{>0}$), and the quotient is isomorphic to $\GL_{2}
(\C)/A\cdot\U (2) \simeq \fH_{3}$, where $A$ denotes the diagonal
subgroup of $\GL_{2} (\C)$.

Our symmetric space $X$ is then built from two copies of $C $,
reflecting the structure of $\iota$. More precisely, let $\V = V\times
V$ and $\cC = C\times C$.  Again $\cC$ is preserved by homotheties,
and we have an diffeomorphism
\begin{equation}\label{eq:diffeo}
\cC / \R_{>0} \longisomto X = G/A_{G}K,
\end{equation}
where $G$, $A_{G}$, $K$ are as in \S\ref{ss:caf}.

\subsection{}
Now we introduce an $F$-structure into the picture.  Let 
$F^{+}\subset F$ be the real quadratic subfield $\Q (\sqrt{5})$.  Then a \emph{binary
Hermitian form over $F$} is a map $\phi\colon F^2 \to F^{+}$ of the form
\[\phi(x,y) = a x \bar{x} + b x \bar{y} + \bar{b}\bar{x} y + c y \bar{y},\]
where $a,c \in F^{+}$ and $b\in F$.  Note that $\pphi = \phi + \phi'$
takes values in $\QQ$.  Indeed, $\pphi$ is precisely the composition
$\Tr_{F^{+}/\QQ} \circ \phi$, and by choosing a $\QQ$-basis for $F$,
$\pphi$ can be viewed as a quaternary quadratic form over $\QQ$.  In
particular, it follows that $\pphi(\OO^2)$ is discrete in $\QQ$.

The \emph{minimum of $\phi$} is \[m(\phi)=\inf_{v \in \OO^2 \setminus
\{ 0\}} \pphi(v).\] A vector $v\in \OO^2$ is \emph{minimal vector} for
$\phi$ if $\phi(v)=m(\phi)$.  The set of minimal vectors for $\phi$ is
denoted $M(\phi)$.  A Hermitian form over $F$ is \emph{perfect} if it
is uniquely determined by $M(\phi)$ and $m(\phi)$.  

\subsection{} We now recall the explicit reduction theory of Koecher
\cite{koecher} and Ash \cite{A} that generalizes work of \Vor on
rational positive-definite quadratic forms \cite{voronoi1}.  Although
these constructions can be done in more generality, we only work with
$\GL_{2}$ over our field $F$. 

Recall that $\V = V\times V$ and $\cC = C\times C$.  Let $q\colon
F^{2}\rightarrow \V$ be the map defined by
\begin{equation}\label{eqn:queue}
q (v) = (vv^{*}, v'v'^{*}).
\end{equation}
Here we view $v$ as a column vector, and $*$ means complex conjugate
transpose.  The restriction of $q$ to $\O^{2}\smallsetminus \{0 \}$
defines a discrete subset $\Xi$ of $\bcC$, the closure of $\cC$ in $\V
$.  Let $\Pi$ be the closed convex hull in $V\times V$ of $\Xi$.  Then
$\Pi$ is an infinite polyhedron known as the \emph{\Vor polyhedron}.
It comes equipped with a natural action of $\GL_{2} (\OO)$.  Modulo
this action $\Pi$ has finitely many faces, and the top-dimensional
faces are in bijection with the perfect quadratic forms over $F$. 

Let $\Sigma$ be the collection of cones on the faces of $\Pi$.  The
set $\Sigma$ forms a \emph{$\Gamma$-admissible polyhedral
decomposition} in the sense of \cite{A}; in particular $\Sigma$
is a fan and admits an action of $\GL_{2} (\OO)$.  When intersected
with the cone $\cC$, the cones in $\Sigma$ provide an explicit
reduction theory for $\GL_{2} (\OO)$ in the following sense.  Any
point $x\in \cC $ is contained in a unique $\sigma (x) \in \Sigma$,
and the set 
\[
\{\gamma \in \GL _{2} (\O ) \mid \gamma \cdot \sigma (x)
= \sigma (x) \}
\]
is finite.  There is also an explicit algorithm to
determine $\sigma (x)$ given $x$, the \emph{\Vor \ reduction
algorithm} \cite{voronoi1, Gmod}.

\subsection{} Every cone $\sigma \in \Sigma$ is preserved by
homotheties, and thus defines a subset in $X$ via \eqref{eq:diffeo}.
We call these subsets \emph{\Vor cells}.  One can think of the \Vor cells
as providing a polytopal tessellation of $X$, although some faces of
the polytopes might be at infinity.  Because of this latter point it
is somewhat awkward to use $\Sigma$ directly to compute cohomology,
although there is a workaround.

According to \cite{A2}, there is a deformation retraction $\cC\to \cC
$ that is equivariant under the actions of both $ \GL_2(\OO)$ and the
homotheties.  Its image modulo homotheties is the \emph{well-rounded
retract} $W$ in~$X$.  The well-rounded retract is contractible, and we
have $H^{*}(\Gamma\bs X; \C) \simeq H^{*} (\Gamma \bs W; \C)$.
Moreover, the quotient $\Gamma\bs W$ is compact.

The well-rounded retract~$W$ is naturally a locally finite cell
complex.  The group $\GL_{2} (\OO )$ preserves the cell structure, and
the stabilizer of each cell in $\GL _{2} (\OO )$ is finite.  One can
show that the cells in $W$ are in a one-to-one, inclusion-reversing
correspondence with the cones in the \Vor fan $\Sigma $ and thus with
the \Vor cells.  This makes it possible to use either the cells in $W$
or the \Vor cells to compute cohomology.  Section 3 of
\cite{computation} gives a very detailed description of how to use $W$
to compute $H^{*} (\Gamma \bs W; \C)$.\footnote{More precisely, in
\cite[\S3]{computation} the authors work with the equivariant
cohomology $H_{\Gamma}^{*} (W; \C)$, but this is isomorphic to $H^{*}
(\Gamma \bs W; \C)$ since $\C$ has characteristic zero.}

\subsection{}
The structure of $\Pi$ in our case has been explicitly determined by
one of us (DY) \cite{Yascyclotomic}. 

Modulo the action of $\GL_{2} (\OO)$, there is one perfect form $\phi$,
represented by the matrix
\[A_\phi=\frac{1}{5}\mat{
\z^3 + \z^2 + 3 & \z^3 - \z^2 + \z - 1 \\
-2\z^3 - \z - 2 & \z^3 + \z^2 + 3 }.\] 

The perfect form $\phi$ has $240$ minimal vectors.  It is clear that
if $v \in M(\phi)$ then $\tau v \in M(\phi)$ for any torsion unit
$\tau \in \OO$; modulo torsion units there are $24$ minimal vectors.
Let $\omega$ denote the unit $\z+\z^2$.  Then modulo torsion
the minimal vectors for $\phi$ are
\begin{multline}\label{eq:perfectmin} 
\vect{
            -\z + 1 \\ \z^3 + 1
            },
            \vect{
            -\z^3+1 \\ 1
            },
            \vect{
            1 \\ -\omega
            },
            \vect{
            1 \\ -\z^2
            },
            \vect{
            1 \\ 0
            },
            \vect{
            1 \\ \z^3
            },
            \vect{
            1 \\ -\z^2 + 1
            },
            \vect{
            1 \\ 1
            },
            \vect{
            1 \\ \z^3 + 1
            },\\
            \vect{
            1 \\ \z + 1
            },
            \vect{
            1 \\ \z^3 + \z + 1
            },
            \vect{
            1 \\ -\z^4
            },
            \vect{
            \omega^{-1} \\ \z^4
            },
            \vect{
            \omega^{-1} \\ \z^4-1
            },
            \vect{
            \omega^{-1} \\ -1
            },
            \vect{
            \omega^{-1} \\ -\z^3 - 1
            },
\\
            \vect{
            \omega^{-1} \\ -\z^3 - \z^2 - 1
            },
            \vect{
            \omega \\ \omega + 1
            },
            \vect{
            \omega \\ -\z^3
            },
            \vect{
            \omega \\ 0
            },
            \vect{
            \omega \\ \z^2
            },
            \vect{
            \omega \\ \omega
            },
            \vect{
            0 \\ 1
            },
            \vect{
            0 \\ \omega
            }
. 
\end{multline}
In $\bcC $ these become $24$ points defining an $8$-cone, which
represents the unique top-dimensional cone in $\Sigma$ modulo $\GL_{2}
(\O)$.  Moreover, one can compute the rest of the cones in $\Sigma$ 
modulo $\GL_{2} (\OO)$.  One finds $5$ $\GL_2(\OO)$-classes of
$7$-cones, $10$ classes of $6$-cones, $11$ classes of $5$-cones, $9$
classes of $4$-cones, $4$ classes of $3$-cones, and $2$ classes of
$2$-cones.  We refer to \cite{Yascyclotomic} for details.

\section{Hecke operators}\label{s:ho}

\subsection{}
The \Vor fan $\Sigma$ gives us a convenient model to compute
cohomology, but unfortunately one cannot use it directly to compute
the action of the Hecke operators.  The problem is that the Hecke
operators, when thought of as Hecke correspondences acting
geometrically on the locally symmetric space $\Gamma \backslash X$, do
not preserve the tessellation corresponding to $\Sigma$.  To address
this problem, we introduce another complex computing the cohomology,
the \emph{sharbly complex} $S_{*}$ \cite{ash.sharb}.
 
Given any nonzero $v\in F^{2}$, let $R (v)\subset \V$ be the ray
through the point $q (v)$ from \eqref{eqn:queue}.  Let $S_k$, $k\geq 0$, be
the $\Gamma$-module $A_{k}/C_{k}$, where (i) $A_{k}$ is the set of formal
$\ZZ$-linear sums of symbols $\vv=[v_1, \cdots, v_{k+2}]$, (ii) each
$v_i$ is a nonzero element of $F^2$, and (iii) $C_{k}$ is the submodule
generated by
\begin{enumerate}
\item $[v_{\sigma(1)}, \cdots, v_{\sigma(k+2)}]-\sgn(\sigma)[v_1,
\cdots, v_{k+2}]$, where $\sigma$ is a permutation on $k+2$ letters, 
\item \label{rayrelation} $[v, v_2, \cdots, v_{k+2}] - [w, v_2, \cdots v_{k+2}]$ if $R (v)
=  R(w)$, and 
\item $[v]$, if $v$ is \emph{degenerate}, i.e., if $v_1, \cdots ,
v_{k+2}$ are contained in a hyperplane.
  
\end{enumerate}  
We define a boundary map $\partial\colon S_{k+1} \to S_{k}$ by 
\begin{equation}\label{eq:boundary}
\del [v_1, \cdots , v_{k+2}] =\sum_{i=1}^{k+2}(-1)^i[v_1, \cdots, \hat{v}_i,\cdots , v_{k+2}].
\end{equation}
This makes $S_{*}$ into a complex.  Note that $S_{*}$ is indexed as a
homological complex, i.e.~the boundary map has degree $(-1)$. We
remark that the our definition is slightly different from that of
\cite{ash.sharb}.  In particular the complex in \cite{ash.sharb} uses
unimodular vectors over $\OO^{2}$ and does not include the relation
(\ref{rayrelation}).  However it is easy to see that the complexes are
quasi-isomorphic.

The basis elements $\vv =[v_1, \cdots, v_{k+2}]$ are called
\emph{$k$-sharblies}.  Our field $F$ has class number $1$, and so
using the relations in $C_k$ one can always find a representative for
$\vv$ with each $v_{i}$ a primitive vector in $\OO^{2}$.  In
particular, one can always arrange that each $q(v_i)$ is a vertex of
$\Pi$.  When such a representative is chosen, the $v_i$ are unique up
to multiplication by a torsion unit in $F$.  In this case the
$v_i$---or by abuse of notation the $q(v_i)$---are called the
\emph{spanning vectors} for $\vv$.  We say a sharbly is
\emph{\Vor{}-reduced} if its spanning vectors are a subset of the
vertices of a \Vor{} cone.

The geometric meaning of this notion is the following.  Each sharbly
$\vv $ with spanning vectors $v_{i}$ determines a closed cone $\sigma
(\vv)$ in $\bcC$, by taking the cone generated by the points $q
(v_{i})$.  Then $\vv$ is \Vor -reduced if and only if $\sigma (\vv)$ is
contained in some \Vor cone.  It is clear that there are finitely many
\Vor-reduced sharblies modulo $\Gamma$.  Not every cone $\sigma (\vv)$
is actually a cone in the fan $\Sigma$, and not every cone in $\Sigma$
has the form $\sigma (\vv)$.  However, as we will see, this causes no
difficulty in our computations.

We can also use the spanning vectors to measure how ``big'' a
$0$-sharbly $\vv$ is: we define the \emph{size of $\vv$}
$\Size (\vv)$  to be the absolute value of the norm determinant of the $2\times 2$
matrix formed by spanning vectors for $\vv$.  By construction $\Size$
takes values in $\ZZ_{>0}$.  In the classical picture for $\bG =
\GL_2/\QQ$, there is only one (up to conjugation and scaling) perfect
form, and it has minimal vectors $e_1, e_2, e_1+e_2$.  Thus for $F =
\QQ$, a $0$-sharbly is \Vor-reduced if and only if it has size $1$,
and a $1$-sharbly is \Vor-reduced if and only if its boundary is
consists of $0$-sharblies of size 1.  For $\GL_{2}$ over general
number fields, the size of a $0$-sharbly $\vv$ is related to whether
or not $\vv$ is \Vor-reduced, but in general there exist \Vor-reduced
$0$-sharblies with size $>1$.

We now consider our field $\QQ(\z)$.  The vertices of a fixed
top-dimensional \Vor cone are given in \eqref{eq:perfectmin}.  Using
this data one can check that a non-degenerate $0$-sharbly is
\Vor-reduced if and only if it has size $1$ or $5$.  For $k > 1$, the
relationship between size and \Vor-reduced $k$-sharblies is more
subtle, but a necessary condition is that each of the sub
$0$-sharblies must have size $1$ or $5$.

The boundary map \eqref{eq:boundary} commutes with the action of
$\Gamma$, and we let $S_*(\Gamma)$ be the homological complex of
coinvariants.  Note that $S_*(\Gamma )$ is infinitely generated as a
$\ZZ \Gamma$-module.  One can show, using Borel--Serre duality
\cite{BS}, that
\begin{equation}\label{eq:iso}
H_{k} ((S_{*}\otimes \CC) (\Gamma))  \longisomto H^{6-k} (\Gamma ; \CC
)
\end{equation}
(cf.~\cite{ash.sharb}).  Moreover, there is a natural action of the
Hecke operators on $S_*(\Gamma)$ (cf.~\cite{experimental}).  We
note that the \Vor-reduced sharblies form a \emph{finitely generated}
subcomplex of $S_{*} (\Gamma)$ that also computes the cohomology of
$\Gamma$ as in \eqref{eq:iso}.  This is our finite model for the
cohomology of $\Gamma$.

\subsection{}
The complex of \Vor-reduced sharblies is not stable under the action 
of Hecke operators.  Thus in order to use the subcomplex of \Vor-reduced 
sharblies to compute Hecke operators, one needs a ``reduction algorithm'' 
for representing the class of a sharbly that is not \Vor-reduced as a sum 
of \Vor-reduced sharblies.  We employ a method analogous to the one described 
in \cite{ants} for real quadratic fields, adapted for the field $F = \QQ(\z)$.

For the convenience of the reader, we recall some of the key points.
By \eqref{eq:iso}, in order to compute cohomology classes in
$H^5(\Gamma;\CC)$, we need to reduce $1$-sharblies.  Specifically, a
cohomology class can be thought of as a linear combination of
\Vor-reduced $1$-sharblies, and the Hecke action sends a \Vor-reduced
$1$-sharbly to a $1$-sharbly that is no longer \Vor-reduced.  The
reduction algorithm an iterative process which proceeds by replacing a
$1$-sharbly that is not \Vor-reduced by a sum of $1$-sharblies that
are closer (in a sense described below) to being \Vor-reduced.

As described above, the boundaries of \Vor-reduced sharblies have
boundaries components of size $0$, $1$, or $5$, and so size gives a
coarse measure of how bad, or far from \Vor-reduced, a $1$-sharbly is.
We describe the reduction of a generic bad $1$-sharbly here; the other
special cases are treated with analogous modifications of \cite{ants}.

\begin{figure}
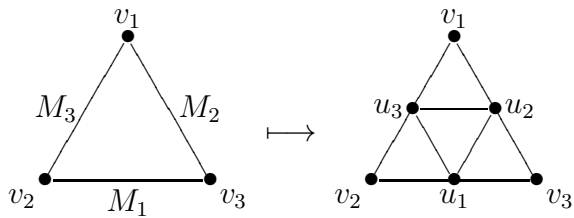
\label{fig:3split}
\[{
\xy /r0.5in/:="a",
{\xypolygon3"A"{~>{-}{\bullet}}},
"A1" *+!D{v_1},"A2" *+!UR{v_2},"A3" *+!UL{v_3},
"A1";"A2" **@{} ?*{}="c1" *+!R{M_3},
"A2";"A3" **@{} ?*{}="c2" *+!U{M_1},
"A3";"A1" **@{} ?*{}="c3" *+!L{M_2},
\endxy} \longmapsto{
\xy /r0.5in/:="a",
{\xypolygon3"A"{~>{-}{\bullet}}},
"A1" *+!D{v_1},"A2" *+!UR{v_2},"A3" *+!UL{v_3},
"A1";"A2" **@{} ?*{\bullet}="c1" *+!R{u_3},
"A2";"A3" **@{} ?*{\bullet}="c2" *+!U{u_1},
"A3";"A1" **@{} ?*{\bullet}="c3" *+!L{u_2},
"c1";"c2" **@{-},"c2";"c3" **@{-},"c3";"c1" **@{-}
\endxy}\]
\caption{Reduction of generic $1$-sharbly.}
\end{figure}

A generic $1$-sharbly $\vv$ that is not \Vor-reduced has boundary components 
that have large size, and so can be thought of as a triangle 
$ \vv = [v_1,v_2,v_3]$ such that $\Size([v_i,v_j]) \gg 0$.  We split each edge 
by choosing \emph{reducing points}
$u_1,u_2$, and $u_3$ and forming three additional edges
$[u_1,u_2],[u_2,u_3]$, and $[u_3,u_1]$.  We then replace $T$ by the
four $1$-sharblies
\begin{equation}\label{eq:generic1sharbly}
[v_1,v_2,v_3]\longmapsto[v_1,u_3,u_2]+[u_3,v_2,u_1]+[u_2,u_1,v_3]+[u_1,u_2,u_3]
\end{equation}
as seen in Figure~\ref{fig:3split}.  Choosing the reducing points uses the 
\Vor polyhedron.  Specifically, the spanning vectors of the $0$-sharbly 
$[v_i,v_j]$ are points in the $8$-dimensional vector space $\V$.  The
 barycenter $b$ of the line joining these points lies in a \Vor cone $\sigma$. 
 The cone $\sigma$ lies between the cone containing $v_i$ and the cone 
containing $v_j$, and so the vertices of $\sigma$ form the candidates for 
reducing points for the $0$-sharbly $[v_i,v_j]$.  We choose the reducing point 
$u$ so that the sum $\Size([v_i,u]) + \Size([u,v_j])$ is minimized.  Note that 
we have not proved that this process decreases size, so we are not guaranteed 
that 
\begin{equation}\label{eq:sizes}
\Size([v_i,v_j]) >\max(\Size([v_i,u]),\Size([u,v_j])).
\end{equation}
Nor are we guaranteed that the $1$-sharbly at the far right of
\eqref{eq:generic1sharbly} is closer to being \Vor-reduced than the
original $1$-sharbly $[v_{1}, v_{2}, v_{3}]$.  However, in practice we
find that both of these problems do not arise.  Indeed, the sizes of
the $0$-sharblies on the right of \eqref{eq:sizes} are typically much
smaller than the size of $[v_{i}, v_{j}]$, and the $1$-sharbly
$[u_{1},u_{2},u_{3}]$ is usually quite close to being \Vor-reduced.

Eventually one produces a $1$-sharbly cycle with all edges
\Vor-reduced.  Unfortunately this is not enough the guarantee that the
cycle \emph{itself} is \Vor-reduced.  This situation does not occur
when one works over $\Q$ as in \cite{computation}, and reflects the
presence of units of infinite order.  Some additional reduction steps
are needed to deal with this problem.  The technique is very similar
to reduction step (IV) in \cite[\S3.5]{ants}.

\section{Elliptic curves}\label{s:ec}

\renewcommand{\o}{\mathcal{O}}

\subsection{} In this section we describe how we constructed the
elliptic curve table at the end of the paper (Table \ref{tab:ec}).
The method itself is the most naive and straightforward one can
imagine.  Recall that $\o = \Z[\zeta]$ and for a positive integer $B$,
let
$$S_B = \bigr\{ c_0 + c_1 \zeta + c_2
\zeta^2 + c_3 \zeta^3 \bigm| |c_i| \leq B, 0\leq i \leq 3 \bigr\},$$ be a boxed
grid of size $16B^4$ inside the lattice of algebraic integers in $F$,
centered at the origin.  For each positive integer $N$, there exists a
bound $B=B(N)$ such that every elliptic curve $E$ over $F$ of
conductor having norm at most $N$ has a Weierstrass model $$E: y^2 +
a_1 xy + a_3 = x^3 + a_2 x^{2} + a_4 x + a_6, \qquad \text{with} \qquad
a_1, a_2, a_3, a_4, a_6 \in S_{B(N)}.$$
 
With $N=10^4$, for example, having found all modular forms which
should correspond to elliptic curves having conductor of norm at most
$N$, we could in principle produce a proof that no elliptic curves not
predicted to exist from the cohomology data up to that level exist as
well as finding the predicted curves.  The bound $B(N)$, however, is
so large as to make this not a practical exercise at the
moment.  Under the assumption of certain conjectures (the ABC
conjecture, for example), one can obtain a much smaller conditional
bound $B^*(N)$, but even this would be far too large to carry out the
proof.

The question we posed for ourselves, therefore, was much a more
practical one: (i) can we perform a reasonable search that finds an
elliptic curve of the predicted conductor matching each rational
cuspidal eigenclass that was found, and (ii) can we in the process
show that though the search is not exhaustive, no unpredicted elliptic
curves appear?

\subsection{} We therefore sifted through curves whose coefficients
$a_i$ lie in the box $S_1$, keeping only those whose discriminants
have modestly sized norm, then filtering those remaining for having
conductor of small norm.  Every curve that was found with conductor
having norm less than $10^4$, matched up with a rational cuspidal
eigenclass from Tables \ref{tab:bettis} and
\ref{tab:qclasses}.  For each of these curves, we then computed
the coefficients $a_\fq = \Norm(\fq) + 1 - | E(\F_\fq)|$ and found
that these matched the Fourier coefficients of the corresponding form
for as many $\fq$ as were computed on the cohomology side. Of course
the computation of $a_\fq$ on the elliptic curve side is very rapid,
so in this way we are able to produce predicted Frobenius eigenvalues
of the modular forms for quite large primes.

\subsection{}
 
Before proceeding with the box search above---which in the end was the
most effective method we could find---we applied another technique.
Although this technique is less systematic, it provides a strategy for
answering a slightly different question: Suppose $F$ is a number field
with class number 1, and an oracle predicts the existence of an
elliptic curve over $F$ of a certain conductor $\mathfrak{n}$, where
$\fn$ is a square-free ideal of fairly small norm, assumed to be odd
for simplicity.  What are some ways in which one can attempt to find a
Weierstrass model for this putative curve and thereby confirm the
prediction of the oracle?  Though we did not need to use it to find
the curves we needed, we discuss a method for answering this question
in case it may serve in another context.

The idea is to use the Frey--Hellegouarch construction of elliptic
curves.  Namely, suppose $u+v=w$ is an equation of $S$-units where $S$
is the set of primes of $\o_F$ dividing $\mathfrak{n}$; this means
that $u,v,w \in \o_F$ are not divisible by primes outside $S$.  Let
$S'$ be the union of $S$ with the set of $\o_F$-primes dividing $2$.
Then, the curve $E_{u,v}$ given by the model $y^2 = x ( x- u) (x+v)$
has good reduction away from $S'$.  By imposing congruences on $u$ and
$v$, we may guarantee that $E_{u,v}$ has good reduction outside $S$.
To speed up the computation, we note that if $\xi$ is an $S$-unit,
then $E_{\xi u, \xi v}$ is a quadratic twist by $\sqrt{\xi}$ of
$E_{u,v}$.  It's also easily seen that curves obtained from
re-orderings such as $E_{v,u}$, $E_{u,-w}$ etc. are also at most
quadratic twists by $\sqrt{-1}$.  Thus, it's convenient to search over
curves $E_{1,\varepsilon}$ and its quadratics twists by square roots
of $S$-units.  In the case of $F=\Q(\zeta)$,
$\o_F^\times/{\o_F^\times}^2$ is generated by $\langle -1, 1+\zeta
\rangle$, so for each $S$-unit equation $1+\varepsilon = \rho$, we get
four curves $E_{1,\varepsilon}$, $E_{-1,-\varepsilon}$, $E_{1+\zeta,
(1+\zeta)\varepsilon}$, and $E_{-(1+\zeta),-(1+\zeta)\varepsilon}$.

One can experimentally search for $S$-unit equations via a similar
grid search as above.  Namely, one finds a basis $\xi_1, \ldots,
\xi_r$ for the group of $S$-units and for a given bound $B$ searches
over integer $r$-tuples $(m_i)_{i=1}^r$ satisfying $|m_i| \leq B$ to
see if $\varepsilon = \prod \xi_i^{m_i}$ yields an $S$-unit equation
$1+\varepsilon=\rho$ by first filtering out those with unsuitable
$\Norm(1+\varepsilon)$.  For example, if $\mathfrak{n} = (\nu)$ is a
principal prime ideal, we can take a basis $\xi_1, \ldots \xi_r$ of
$\o_F^\times$ and check whether $\Norm(\nu - \prod_i \xi_i^{m_i})=\pm
1$.
 
As an example over $F=\Q(\zeta)$, the unit group is generated by
$-\zeta, 1+\zeta$.  If we take $1+\varepsilon=\rho$ where
$\varepsilon=-\zeta^3(1+\zeta)^{24}$, then the curve
$E_{1,\varepsilon}$ has discriminant $\Delta = 2^{12}3^4 5 \zeta^4
(1+\zeta)^{72}$, conductor $\mathfrak{n}= (3-3\zeta)$ of norm $405$,
which is the second conductor listed in Table \ref{tab:levels}.

As another example, the 5-unit equation $u + v = w$ where
$u=\zeta^3(1+\zeta)^{-1}=\zeta+\zeta^{-1}$ and $v= \zeta^2(1+\zeta)$
is especially nice because $v-u=1, -uv=-1$ yields that the curve
$E_{u,v}$ is $y^2=x^3 + x^2 - x$ which descends to $\Q$.  It has
conductor of norm $1280$, the sixth conductor listed in Table
\ref{tab:levels}.

\section{Results}\label{s:data}

\newcommand{\ourprime}{12379}

\subsection{} In this section we present our computational data, both
on the cohomology and elliptic curve sides.  Our programs were
implemented in Magma \cite{magma}.  We remark that in the cohomology
computations, following a standard practice (cf.~\cite{computation})
we did not work over the complex numbers $\C$, but instead computed
cohomology with coefficients in a large finite field $\F_{\ourprime}$.
This technique was used to avoid the precision problems in
floating-point arithmetic.  Since we do not expect $\Gamma_{0} (\fn)$
to have $\ourprime$-torsion, we expect that the Betti numbers we
report coincide with those one would compute for the group cohomology
with $\C$-coefficients.  As a check, we reran some computations with
coefficients in finite fields over other large primes, and found the
same Betti numbers each time.  Thus we believe we are actually
reporting the dimensions of $H^{5} (\Gamma_{0} (\fn)\backslash X; \C)$.

\subsection{Cuspidal cohomology}
Our first task was to identify those levels with nonzero cuspidal
cohomology.  We first experimentally determined the dimensions of the
subspace $H^{5}_{\Eis}$ spanned by \emph{Eisenstein cohomology
classes} \cite{harder.gl2}.  Such classes are closely related to
Eisenstein series.  In particular the eigenvalue of $T_{\fq}$ on these
classes equals $\Norm (\fq) +1$.  We expect that for a given level
$\fn$, the dimension of the Eisenstein cohomology space depends only
the factorization type of $\fn$.  Thus initially we used some Hecke
operators applied to cohomology spaces of small level norm to compute
the expected Eisenstein dimension for small levels with different
factorization types.  The result can be found in
Table~\ref{tab:eisenstein}.

After compiling Table~\ref{tab:eisenstein}, we computed cohomology for
a larger range of levels and looked for Betti numbers in excess of
that in Table~\ref{tab:eisenstein}.  We were able to compute $H^{5}$
for all levels $\fn$ with $\Norm(\fn) \leq 4941$.  For $\fn =\fp$
prime we were able to carry the computations further to $\Norm (\fp)
\leq 7921$.  Table~\ref{tab:levels} shows the norms of the levels
$\fn$ with nonzero cuspidal cohomology and generators of $\fn$ we
used.  It turns out that modulo the action of Galois each cuspidal
space can be uniquely identified by the norm of the level, except when
$\Norm (\fn) = 3641$.  In this case there are two levels up to Galois
with nonzero cuspidal cohomology; we call them 3641a and 3641b.  The
dimensions of the cuspidal subspaces $H^{5}_{\cusp}$ are given in
Table~\ref{tab:bettis}.

\subsection{Hecke operators} Next we computed the Hecke operators and
looked for eigenclasses with rational eigenvalues.  These computations
were quite intensive.  For all levels we were able to compute at least
up to $T_{\fq}$ with $\fq \subset \O$ prime satisfying $\Norm
(\fq)\leq 41$; at some smaller levels, such as $\Norm (\fn) = 701$, we
computed much further.  At the largest levels ($\Norm (\fn) = 4455,
4681, 6241, 7921$) the computation was so big that our implementation
could not compute any Hecke operators.  Table~\ref{tab:heckereps}
gives our choices of generators for the ideals $\fq$.

For all levels except for one, the cuspidal cohomology split into
$1$-dimensional rational eigenspaces.  We give some eigenvalues for
the rational eigenclasses in Table \ref{tab:qclasses}.  The remaining
level --- norm 3721 --- is $2$-dimensional with Hecke eigenvalues
generating the field $F^{+} = \Q (\sqrt{5})$.  The characteristic
polynomials can be seen in Table \ref{tab:3721-table}.

\subsection{Elliptic curves over $F$} Now we give motivic explanations
for all the cuspidal cohomology we found.

Thirteen of the eigenclasses in Table \ref{tab:qclasses} have the
property that their eigenvalues $a_{\fq}$ differ for at least two
primes $\fq , \fq '$ lying over the same prime in the subfield $F^+$.
Hence we expect these classes to correspond to elliptic curves over
$F$.  Using the techniques described in \S\ref{s:ec}, we were able to
find elliptic curves $E/F$ such that for all primes $\fq$ of good
reduction, the identity $a_{\fq} = \Norm (\fq) + 1 - |E (\F_{\fq})|$
held for every Hecke operator we computed.  Equations for these curves
are given in Table \ref{tab:ec}.

Although we were unable to match the remaining eigenclass, namely
the second labelled 3641b, to an elliptic curve over $F$, a curve
matching this class was found by Mark Watkins (Appendix A).  

\subsection{The remaining eigenclasses} All the other eigenclasses
Tables \ref{tab:qclasses} and \ref{tab:3721-table} can be accounted
for either by elliptic curves over $\Q$, elliptic curves over $F^{+}$,
``old'' cohomology classes coming from lower levels, or other Hilbert
modular forms over $F^{+}$.  We indicate briefly what happens.

\subsubsection{Elliptic curves over $\Q$} The eigenclasses at 400,
405, 1280, 1296, 4096, and one of the eigenclasses at 2025, correspond
to elliptic curves over $\Q$ that can readily be found in Cremona's
tables \cite{cremona}.  In all cases, there are actually \emph{two}
rational elliptic curves that are not isogenous over $\Q$ but produce
the same eigenvalue data when considered as curves over $F$; the
curves in these pairs are quadratic twists by $5$ of each other that
become isomorphic over $F^{+}$. For instance, at 400 the two curves
are 50A1 and 50B3 (in the notation of \cite{cremona}).

\subsubsection{Elliptic curves over $F^{+}$} The eigenclasses at 605,
961, 1681, 1805, 2401, and 4205 correspond to elliptic curves over
$F^{+}$.  The class at 2401 already appears in \cite{dembele}; the
others were verified using software written by Demb\'el\'e.  As an
example, the three eigenclasses at 4205 correspond to three cuspidal
parallel weight 2 Hilbert modular newforms of level $\fp_{5}\fp_{29}
\subset \O_{F^{+}}$. Although we were unable to compute Hecke
operators at 6241 and 7921, we expect that these classes correspond to
elliptic curves given in \cite{dembele}.

\subsubsection{Old classes} There are two-dimensional eigenspaces at
2000, 2025, 3025, 3505, 4400, and 4455 on which the Hecke operators we
computed act by scalars.  These subspaces correspond to curves
appearing at lower levels.  For example, the classes at 2000 and 4400
correspond to the classes that already appeared at 400.  We note
that 2000, 2025, 4400, and 4455 correspond to elliptic curves over
$\Q$, while 3025 corresponds to an elliptic curve over $F^{+}$ (seen
in Table \ref{tab:qclasses} at 605) and 3505 to a curve over $F$ (seen
in Table \ref{tab:qclasses} at 701).

\subsubsection{Other Hilbert modular forms} There are two eigenclasses
remaining, namely the class at 3721 with eigenvalues in $F^{+}$ and
the third eigenclass $\xi$ at 3025 with eigenvalues in $\Q$.  Both can be
attributed to Hilbert modular forms of parallel weight 2 attached to
abelian surfaces.

For 3721, the characteristic polynomials match those of a parallel
weight 2 Hilbert modular newform of level $\fp_{61}\subset
\O_{F^{+}}$.

The class $\xi$ at 3025 is perhaps the most interesting of all, other than
the classes matching elliptic curves over $F$, since it gives an example of
a \emph{fake elliptic curve} in the sense of \cite{cremona.fake}.  Let
$\fm \subset \O_{F^+}$ be the ideal $\fp_{5}^{2}\fp_{11}$.  The space of
parallel weight 2 Hilbert modular newforms of level $\fm$ contains an
eigenform $g$ with Hecke eigenvalues $a_{\fq}$ in the field $F^{+}$.
For any prime $\fq \subset \O_{F^{+}}$, let $q\in \Z $ be the prime
under $\fq$.  Then we have $a_{\fq} (g) =0$ if $q=5$, and 
\begin{equation}\label{eq:fake}
a_{\fq} (g) \in \begin{cases}
\Z&\text{if $q=1\bmod 5$},\\
\Z\cdot\sqrt{5} & \text{if $q=2,3,4 \bmod 5$}.
\end{cases}
\end{equation}
Table~\ref{tab:3025-table} gives some eigenvalues of $g$.  The
conditions \eqref{eq:fake} imply that there is a quadratic character
$\varepsilon$ of $\Gal (F/F^{+})$ such that the $L$-series $L (s,g) L
(s, g\otimes \varepsilon)$ agrees with the $L$-series attached to our
eigenclass $\xi$.  Indeed, following \cite{cremona.fake}, if $\fq \subset
\O_{F^{+}}$ splits in $F$ as $\fr \cdot \bar\fr$ (respectively, remains
inert in $F$), then we should expect the Hecke eigenvalues of $g$ and
$\xi$ to be related by 
\[
a_{\fr} (\xi ) = a_{\bar\fr} (\xi) = a_{\fq} (g) \quad \text{(split)}
\]
and 
\[
a_{\fq} (\xi) = a_{\fq}(g)^{2} - 2\Norm_{F^{+}/\Q} (\fq ) \quad \text{(inert)}.
\]
Comparison of Tables~\ref{tab:qclasses} and~\ref{tab:3025-table}
shows that this holds.

\appendix
\section{Elliptic curves with good reduction outside a given set (Mark
Watkins)}\label{sec:appendix}

The method to find a curve with good reduction outside a finite set is
outlined in Cremona--Lingham \cite{cl}, though much of this was well-known
to experts in prior times. In our specific case, we can make some
additional simplifications and/or modifications.

Since the primes $S=\{\p_{11},\p_{331}\}$ that divide the level are
exactly the same as the primes that divide the discriminant of the
elliptic curve, we immediately have that $\Delta=(-1)^a u^b e_{11}^c
e_{331}^d$ where $u=1+\zeta_5^2+\zeta_5^3$ is a unit,
$\p_{11}=(e_{11})$ and $\p_{331}=(e_{331})$ are principalisations, and
$a\in\{0,1\}$, $0\le b\le 11$, and $c,d\ge 1$.

The formul\ae\ $j=c_4^3/\Delta$ and $j-1728=c_6^2/\Delta$ then give us
various divisibility conditions.  For instance, upon noting the
triviality of the class group of ${\bf Q}(\zeta_5)$, Proposition of
3.2 of \cite{cl} implies that $w=j^2(j-1728)^3$ must have $6|v_\p(w)$
for all primes~$\p$ other than $\p_{11}$ and~$\p_{331}$.  It is a
standard problem in algorithmic number theory to list all possible
such~$w\in{\bf Q}(\zeta_5)$ up to 6th powers, and then for each~$w$ we
are left to find $S$-integral points on the curve $E(w):
Y^2=X^3-1728w$.

We can work more directly in our case, and note that
$j=c_4^3/\Delta=1728+c_6^2/\Delta$ gives an elliptic curve $E(\Delta):
c_6^2=c_4^3-1728\Delta$ in the unknowns~$c_4$ and~$c_6$.  We can note
that two curves with $\Delta$ differing by a 6th power will give
isomorphic $E(\Delta)$, though in making such a passage we may need to
find $S$-integral points on the resulting curves rather than just
integral points. Also, this allows us to restrict to $0\le b\le 5$
without loss.  We are unable to find the full Mordell--Weil group for
most of the $E(\Delta)$ curves in any event, and so completeness is
impractical.

It turns out that $(a,b,c,d)\in\{(1,3,2,1),(1,5,2,1)\}$ will give the
first and second curves corresponding to 3641b.  To find these, we
tried all possibilities for $(a,b)$ with $(c,d)=(1,1)$, and then
$(c,d)=(2,1)$. Thus we had to try to find the Mordell--Weil group for
24 different elliptic curves (we were succesful for only~7).  The
curves~$E(\Delta)$ all have a 3-isogeny, but this does not seem to be
of much use.

We used the Magma package of Nils Bruin to try to search for points on
the~$E(\Delta)$.  We can get an upper bound on the rank using {\tt
TwoSelmerGroup}, though this is not strictly necessary.  We then
search for points on the elliptic curves using the function {\tt
PseudoMordellWeilGroup} with a {\tt SearchBound} of 100.  This took
about 5 minutes per curve (the search bound is about the 4th power of
this in terms of the norm, as the field is quartic).  In
Table~\ref{table-app} we list the data for the upper bound on the rank
and the number of generators found.

It is natural that we can find more points when the rank is large, as
the points are more likely to be of smaller height. Once we have some
linearly independent points in the Mordell--Weil group, we can find all
integral points that they generate.  Again a provable version of this
is rather technical, and largely unneeded.  We simply took all linear
combinations with coefficients of size not more than 5. This then
gives a set of integral points $(X,Y)$ on~$E(\Delta)$, and from each
we can obtain an elliptic curve with the correct $j$-invariant via
$$j=X^3/\Delta\quad\text{and}\quad
E_\Delta(X):y^2=x^3-{3j\over j-1728}x-{2j\over j-1728}.$$

We can then try to twist away ramification at places outside~$\p_{11}$
and~$\p_{331}$. However, we can also perform a preliminary check on the
traces of Frobenius of the curves $E_\Delta(X)$, as they must match those
from the Hecke operators up to sign if the twisting is to be successful.

We are fortunate in the end, since even though Table \ref{table-app}
contains many missing Mordell--Weil groups, we are still able to find
the two desired curves.  In Table \ref{table-app}, the $(a,b,c,d)$ column gives
the choice of these parameters in the discriminant, the $s$-column
gives the upper bound on the rank from {\tt TwoSelmerGroup}, the
$g$-column gives the number of generators we found via a search up to
na\"i{}ve height~100, and $I$-column gives the number of integral
points we obtained from these when taking small linear combinations of
the generators.

We conclude by giving the Weierstrass equation for the second curve
labelled 3641b:
\begin{multline}\label{eq:3641bcurve}
y^2 +
(\z^2 + 1) xy + \z^2 
= x^3 +  (-\z^3 + \z^2 + \z + 1) x^{2} \\
+ 
(-\z^3 - 82\z^2 + 52\z - 84 )x + (310\z^3 - 366\z^2 + 418\z - 175).
\end{multline}
\bibliographystyle{amsplain_initials}
\bibliography{gunnells-hajir-yasaki}

\providecommand{\bysame}{\leavevmode\hbox to3em{\hrulefill}\thinspace}
\providecommand{\MR}{\relax\ifhmode\unskip\space\fi MR }
\providecommand{\MRhref}[2]{%
  \href{http://www.ams.org/mathscinet-getitem?mr=#1}{#2}
}
\providecommand{\href}[2]{#2}
\begin{thebibliography}{10}

\bibitem{A}
A.~Ash, \emph{Deformation retracts with lowest possible dimension of arithmetic
  quotients of self-adjoint homogeneous cones}, Math. Ann. \textbf{225} (1977),
  no.~1, 69--76.

\bibitem{A2}
\bysame, \emph{Small-dimensional classifying spaces for arithmetic subgroups of
  general linear groups}, Duke Math. J. \textbf{51} (1984), no.~2, 459--468.

\bibitem{ash.sharb}
\bysame, \emph{Unstable cohomology of ${\SL}(n,\mathscr{O})$}, J. Algebra
  \textbf{167} (1994), no.~2, 330--342.

\bibitem{computation}
A.~Ash, P.~E. Gunnells, and M.~McConnell, \emph{Cohomology of congruence
  subgroups of {${\rm SL}\sb 4(\Bbb Z)$}}, J. Number Theory \textbf{94} (2002),
  no.~1, 181--212.

\bibitem{BS}
A.~Borel and J.-P. Serre, \emph{Corners and arithmetic groups}, Comment. Math.
  Helv. \textbf{48} (1973), 436--491, Avec un appendice: Arrondissement des
  vari\'et\'es \`a coins, par A. Douady et L. H\'erault.

\bibitem{bw}
A.~Borel and N.~Wallach, \emph{Continuous cohomology, discrete subgroups, and
  representations of reductive groups}, second ed., Mathematical Surveys and
  Monographs, vol.~67, American Mathematical Society, Providence, RI, 2000.

\bibitem{borel.survey}
A.~Borel, \emph{Introduction to the cohomology of arithmetic groups}, Lie
  groups and automorphic forms, AMS/IP Stud. Adv. Math., vol.~37, Amer. Math.
  Soc., Providence, RI, 2006, pp.~51--86.

\bibitem{magma}
W.~Bosma, J.~Cannon, and C.~Playoust, \emph{The {M}agma algebra system. {I}.
  {T}he user language}, J. Symbolic Comput. \textbf{24} (1997), no.~3-4,
  235--265, Computational algebra and number theory (London, 1993).

\bibitem{bygott}
J.~Bygott, \emph{Modular forms and modular symbols over imaginary quadratic
  fields}, Ph.D. thesis, Exeter, 1999.

\bibitem{cremona2}
J.~E. Cremona, \emph{Hyperbolic tessellations, modular symbols, and elliptic
  curves over complex quadratic fields}, Compositio Math. \textbf{51} (1984),
  no.~3, 275--324.

\bibitem{cremona.fake}
\bysame, \emph{Abelian varieties with extra twist, cusp forms, and elliptic
  curves over imaginary quadratic fields}, J. London Math. Soc. (2) \textbf{45}
  (1992), no.~3, 404--416.

\bibitem{cremona}
\bysame, \emph{The elliptic curve database for conductors to 130000},
  Algorithmic number theory, Lecture Notes in Comput. Sci., vol. 4076,
  Springer, Berlin, 2006, pp.~11--29.

\bibitem{cl}
J.~E. Cremona and M.~P. Lingham, \emph{Finding all elliptic curves with good
  reduction outside a given set of primes}, Experiment. Math. \textbf{16}
  (2007), no.~3, 303--312.

\bibitem{crem.whitley}
J.~E. Cremona and E.~Whitley, \emph{Periods of cusp forms and elliptic curves
  over imaginary quadratic fields}, Math. Comp. \textbf{62} (1994), no.~205,
  407--429.

\bibitem{dembele}
L.~Demb{\'e}l{\'e}, \emph{Explicit computations of {H}ilbert modular forms on
  {${\Bbb Q}(\sqrt{5})$}}, Experiment. Math. \textbf{14} (2005), no.~4,
  457--466.

\bibitem{Fra}
J.~Franke, \emph{Harmonic analysis in weighted {$L\sb 2$}-spaces}, Ann. Sci.
  \'Ecole Norm. Sup. (4) \textbf{31} (1998), no.~2, 181--279.

\bibitem{Gmod}
P.~E. Gunnells, \emph{Modular symbols for {${\bf Q}$}-rank one groups and
  {V}orono\u\i\ reduction}, J. Number Theory \textbf{75} (1999), no.~2,
  198--219.

\bibitem{experimental}
\bysame, \emph{Computing {H}ecke eigenvalues below the cohomological
  dimension}, Experiment. Math. \textbf{9} (2000), no.~3, 351--367.

\bibitem{ants}
P.~E. Gunnells and D.~Yasaki, \emph{Hecke operators and {H}ilbert modular
  forms}, Algorithmic number theory, Lecture Notes in Comput. Sci., vol. 5011,
  Springer, Berlin, 2008, pp.~387--401.

\bibitem{harder.book}
G.~Harder, \emph{Cohomology of arithmetic groups}, book in preparation
  available from Harder's website.

\bibitem{harder.gl2}
\bysame, \emph{Eisenstein cohomology of arithmetic groups. {T}he case {${\rm
  GL}_2$}}, Invent. Math. \textbf{89} (1987), no.~1, 37--118.

\bibitem{koecher}
M.~Koecher, \emph{Beitr\"age zu einer {R}eduktionstheorie in
  {P}ositivit\"atsbereichen. {I}}, Math. Ann. \textbf{141} (1960), 384--432.

\bibitem{li.schwermer.survey}
J.-S. Li and J.~Schwermer, \emph{Automorphic representations and cohomology of
  arithmetic groups}, Challenges for the 21st century ({S}ingapore, 2000),
  World Sci. Publ., River Edge, NJ, 2001, pp.~102--137.

\bibitem{lingham}
M.~Lingham, \emph{Modular forms and elliptic curves over imaginary quadratic
  fields}, Ph.D. thesis, Nottingham, 2005.

\bibitem{sw}
J.~Socrates and D.~Whitehouse, \emph{Unramified {H}ilbert modular forms, with
  examples relating to elliptic curves}, Pacific J. Math. \textbf{219} (2005),
  no.~2, 333--364.

\bibitem{voronoi1}
G.~Vorono\v\i, \emph{Sur quelques propri\'et\'es des formes quadratiques
  positives parfaites}, J. Reine Angew. Math. \textbf{133} (1908), 97--178.

\bibitem{Yascyclotomic}
D.~Yasaki, \emph{Binary {H}ermitian forms over a cyclotomic field}, J. Algebra
  \textbf{322} (2009), 4132--4142.

\end{thebibliography}

\newpage

\begin{table} \begin{tabular}{|l|ccccccccccc|}
\hline%
Factorization of $\fn$& $\fp$  & $\fp^2$&$\fp^3$ & $\fp^4$ &$\fp^5$& $\fp \fq$ & $\fp^2 \fq$ & 
$\fp^3 \fq$ &$\fp^2\fq^2$& $\fp\fq\fr$ &$\fp^2\fq\fr$\\ 
\hline
$\dim H^{5}_{\Eis} (\Gamma_{0} (\fn))$& 3&5&7&9&11&7&11&15&17&15&23\\
\hline
\end{tabular}
\medskip
\caption{\label{tab:eisenstein}
Expected dimension of Eisenstein cohomology $H^{5}_{\Eis} (\Gamma_{0}
(\fn))$ in terms of the prime factorization of $\fn$.  Prime ideals
are denoted by $\fp$, $\fq$, $\fr$.}
\end{table}

\begin{table}[htb]
\begin{center}
\begin{tabular}{|l|r|l|r|l|r|}
\hline
N($\fn$) & generator of $\fn$&N($\fn$) & generator of $\fn$&N($\fn$) & generator of $\fn$\cr
\hline
400 & $2\z^2-4\z+2$& 405 & $-3\z^3-3\z^2-3\z-6$& 605 & $-\z^3+4\z^2-4\z+1$\cr
701 & $-2\z^3-\z^2-3\z-6$& 961 & $-2\z^3-2\z^2+5$& 1280 & $-4\z+4$\cr
1296 & $6$& 1681 & $-6\z^3-7\z^2-7\z-6$& 1805 & $-3\z^3-4\z^2-5\z-8$\cr
2000 & $-2\z^3+6\z^2-6\z+2$& 2025 & $-3\z^3-3\z^2-9$& 2201 & $-2\z^3-6\z^2-7\z-8$\cr
2351 & $-2\z^3-6\z^2-\z-9$& 2401 & $7$& 3025 & $-6\z^3+7\z^2-6\z$\cr
3061 & $-6\z^3-7\z^2-5\z-10$& 3355 & $5\z^3-5\z^2+2\z+3$& 3505 & $-2\z^3-8\z^2+\z-11$\cr
3571 & $-4\z^2-6\z-9$& 3641a & $-2\z^3-5\z^2+4\z-10$& 3641b & $-\z^3+7\z^2-4\z+1$\cr
3721 & $7\z^3+7\z^2+3$& 4096 & $8$& 4205 & $-4\z^3-5\z^2-6\z-10$\cr
4400 & $4\z^3+10\z^2-2\z+8$& 4455 & $-6\z^2-9$& 4681 & $\z^3-8\z^2-7$\cr
5081 & $-2\z^3-5\z^2-5\z-10$& 5101 & $-6\z^3-2\z^2-11$& 6241 & $3\z^3+11\z^2+3\z$\cr
6961 & $-8\z^3-6\z^2-5\z-14$& 7921 & $-11\z^3-\z^2-\z-11$& &\cr
\hline
\end{tabular}
\end{center}
\medskip
\caption{Levels $\fn$ with nontrivial cuspidal cohomology.  Only one representative of each level up to Galois is given.\label{tab:levels}}
\end{table}

\begin{table}[htb]
\begin{center}
\begin{tabular}{|l|c|l|c|l|c|}
\hline
N($\fn$) & dimension&N($\fn$) & dimension&N($\fn$) & dimension\cr
\hline
400 & $1$& 405 & $1$& 605 & $1$\cr
701 & $1$& 961 & $1$& 1280 & $1$\cr
1296 & $1$& 1681 & $1$& 1805 & $1$\cr
2000 & $2$& 2025 & $3$& 2201 & $1$\cr
2351 & $1$& 2401 & $1$& 3025 & $3$\cr
3061 & $1$& 3355 & $1$& 3505 & $2$\cr
3571 & $1$& 3641a & $1$& 3641b & $2$\cr
3721 & $2$& 4096 & $1$& 4205 & $3$\cr
4400 & $2$& 4455 & $2$& 4681 & $1$\cr
5081 & $1$& 5101 & $1$& 6241 & $1$\cr
6961 & $1$& 7921 & $1$& &\cr
\hline
\end{tabular}
\end{center}
\medskip
\caption{Dimensions of cuspidal cohomology.\label{tab:bettis}}
\end{table}

\begin{table}
\begin{center}
\begin{tabular}{|l|r|}
\hline
Prime under $\fq$ & Generator of $\fq$\cr
\hline
 2 & $2$ \cr
\hline
 5 & $-\z+1$ \cr
\hline
 11 & $-\z^3-\z+1$ \cr
 & $\z^3-\z+1$ \cr
 & $-\z^2+\z+1$ \cr
 & $2\z^2+\z+1$ \cr
\hline
 31 & $-2\z+1$ \cr
 & $-2\z^2+1$ \cr
 & $-2\z^3+1$ \cr
 & $2\z^3+2\z^2+2\z+3$ \cr
\hline
 41 & $-\z^3-3\z^2-\z-2$ \cr
 & $-\z^3-\z^2-2\z-3$ \cr
 & $-\z^3-2\z^2-2\z-3$ \cr
 & $2\z^3+3\z^2+\z+2$ \cr
\hline
\end{tabular}
\end{center}
\medskip
    \caption{Choice of primes $\fq$ for Hecke operators\label{tab:heckereps}.}
\end{table}

\begin{table}
\begin{center}\scriptsize
\begin{tabular}{|l|rrrrrrrrrrrrrr|}
\hline
N($\fn$) & 2 & 5 & 11 & 11 & 11 & 11 & 31 & 31 & 31 & 31 & 41 & 41 & 41 & 41\cr
\hline
400 & $\bullet$& $\bullet$& $-3$& $-3$& $-3$& $-3$& $2$& $2$& $2$& $2$& $-3$& $-3$& $-3$& $-3$\cr
405 & $1$& $\bullet$& $-4$& $-4$& $-4$& $-4$& $0$& $0$& $0$& $0$& $10$& $10$& $10$& $10$\cr
605 & $-7$& $\bullet$& $\bullet$& $\bullet$& $\bullet$& $\bullet$& $8$& $-4$& $-4$& $8$& $-6$& $6$& $6$& $-6$\cr
701 & $-1$& $-3$& $3$& $3$& $-6$& $3$& $-4$& $5$& $-4$& $-4$& $6$& $-3$& $6$& $-12$\cr
961 & $1$& $-2$& $4$& $4$& $-4$& $-4$& $\bullet$& $\bullet$& $\bullet$& $\bullet$& $-6$& $-6$& $-6$& $-6$\cr
1280 & $\bullet$& $\bullet$& $0$& $0$& $0$& $0$& $-4$& $-4$& $-4$& $-4$& $6$& $6$& $6$& $6$\cr
1296 & $\bullet$& $-4$& $2$& $2$& $2$& $2$& $-8$& $-8$& $-8$& $-8$& $2$& $2$& $2$& $2$\cr
1681 & $-4$& $-1$& $5$& $5$& $-2$& $-2$& $-10$& $4$& $4$& $-10$& $\bullet$& $\bullet$& $\bullet$& $\bullet$\cr
1805 & $-7$& $\bullet$& $0$& $0$& $0$& $0$& $-4$& $8$& $8$& $-4$& $-6$& $-6$& $-6$& $-6$\cr
2000 & $\bullet$& $\bullet$& $-3$& $-3$& $-3$& $-3$& $2$& $2$& $2$& $2$& $-3$& $-3$& $-3$& $-3$\cr
2000 & $\bullet$& $\bullet$& $-3$& $-3$& $-3$& $-3$& $2$& $2$& $2$& $2$& $-3$& $-3$& $-3$& $-3$\cr
2025 & $1$& $\bullet$& $-4$& $-4$& $-4$& $-4$& $0$& $0$& $0$& $0$& $10$& $10$& $10$& $10$\cr
2025 & $1$& $\bullet$& $-4$& $-4$& $-4$& $-4$& $0$& $0$& $0$& $0$& $10$& $10$& $10$& $10$\cr
2025 & $-8$& $\bullet$& $2$& $2$& $2$& $2$& $-3$& $-3$& $-3$& $-3$& $-8$& $-8$& $-8$& $-8$\cr
2201 & $1$& $-4$& $-3$& $-6$& $-5$& $-2$& $\bullet$& $\bullet$& $\bullet$& $\bullet$& $4$& $-3$& $-6$& $0$\cr
2351 & $3$& $-1$& $-2$& $5$& $-2$& $-2$& $4$& $-3$& $4$& $4$& $0$& $0$& $0$& $7$\cr
2401 & $-8$& $-4$& $-3$& $-3$& $-3$& $-3$& $2$& $2$& $2$& $2$& $2$& $2$& $2$& $2$\cr
3025 & $-7$& $\bullet$& $\bullet$& $\bullet$& $\bullet$& $\bullet$& $8$& $-4$& $-4$& $8$& $-6$& $6$& $6$& $-6$\cr
3025 & $-7$& $\bullet$& $\bullet$& $\bullet$& $\bullet$& $\bullet$& $8$& $-4$& $-4$& $8$& $-6$& $6$& $6$& $-6$\cr
3025 & $-3$& $\bullet$& $\bullet$& $\bullet$& $\bullet$& $\bullet$& $-8$& $2$& $2$& $-8$& $2$& $-8$& $-8$& $2$\cr
3061 & $-3$& $-4$& $-4$& $-3$& $-1$& $-2$& $-2$& $-2$& $-9$& $-6$& $6$& $-4$& $-5$& $3$\cr
3355 & $5$& $\bullet$& $\bullet$& $\bullet$& $\bullet$& $\bullet$& $-4$& $-4$& $-4$& $8$& $-6$& $-6$& $-6$& $-6$\cr
3505 & $-1$& $\bullet$& $3$& $3$& $-6$& $3$& $-4$& $5$& $-4$& $-4$& $6$& $-3$& $6$& $-12$\cr
3505 & $-1$& $\bullet$& $3$& $3$& $-6$& $3$& $-4$& $5$& $-4$& $-4$& $6$& $-3$& $6$& $-12$\cr
3571 & $-5$& $-3$& $-6$& $-2$& $-3$& $$& $-8$& $-5$& $-2$& $0$& $-2$& $8$& $10$& $-3$\cr
3641a & $-7$& $-3$& $\bullet$& $\bullet$& $\bullet$& $\bullet$& $-1$& $-6$& $3$& $11$& $0$& $-2$& $-9$& $-2$\cr
3641b & $-1$& $-3$& $\bullet$& $\bullet$& $\bullet$& $\bullet$& $-1$& $8$& $-7$& $-7$& $-12$& $0$& $9$& $0$\cr
3641b & $7$& $1$& $\bullet$& $\bullet$& $\bullet$& $\bullet$& $7$& $-8$& $-3$& $-3$& $12$& $-8$& $-3$& $-8$\cr
4096 & $\bullet$& $-2$& $-4$& $-4$& $-4$& $-4$& $0$& $0$& $0$& $0$& $2$& $2$& $2$& $2$\cr
4205 & $-4$& $\bullet$& $5$& $5$& $-2$& $-2$& $-10$& $-3$& $-3$& $-10$& $0$& $7$& $7$& $0$\cr
4205 & $-4$& $\bullet$& $-3$& $-3$& $-6$& $-6$& $2$& $5$& $5$& $2$& $0$& $-9$& $-9$& $0$\cr
4205 & $-7$& $\bullet$& $-4$& $-4$& $4$& $4$& $8$& $0$& $0$& $8$& $-6$& $10$& $10$& $-6$\cr
4400 & $\bullet$& $\bullet$& $\bullet$& $\bullet$& $\bullet$& $\bullet$& $2$& $2$& $2$& $2$& $-3$& $-3$& $-3$& $-3$\cr
4400 & $\bullet$& $\bullet$& $\bullet$& $\bullet$& $\bullet$& $\bullet$& $2$& $2$& $2$& $2$& $-3$& $-3$& $-3$& $-3$\cr
5081 & $3$& $-4$& $-4$& $0$& $0$& $-6$& $0$& $0$& $-8$& $0$& $-6$& $6$& $-6$& $-4$\cr
5101 & $-3$& $-3$& $-1$& $0$& $-3$& $-5$& $7$& $-10$& $1$& $-8$& $10$& $4$& $-12$& $-10$\cr
6961 & $1$& $-2$& $0$& $-6$& $-4$& $0$& $-10$& $4$& $-8$& $-2$& $-10$& $-8$& $0$& $-2$\cr
\hline
\end{tabular}
\end{center}
\medskip
\caption{Eigenvalues of cuspidal $\Q$-eigenclasses\label{tab:qclasses}.  For each Hecke operator $T_\fq$ we give the rational prime lying under $\fq$.  The order of the columns corresponds to Table \ref{tab:heckereps}.}
\end{table}

\begin{table}
\begin{center}
\begin{tabular}{|l|r|}
\hline
Prime under $\fq$ & Characteristic polynomial of $T_\fq$\cr
\hline
 2 & $x^2 + 4x - 16$ \cr
\hline
 5 & $x^2 + x - 11$ \cr
\hline
 11 & $x^{2}-20$ \cr
 & $x^{2}-20$ \cr
 & $x^2 + x - 1$ \cr
 & $x^2 + x - 1$ \cr
\hline
 31 & $x^2 - 14x + 44$ \cr
 & $x^2 - 14x + 44$ \cr
 & $x^2 + 3x - 29$ \cr
 & $x^2 + 3x - 29$ \cr
\hline
 41 & $x^2 - 10x + 20$ \cr
 & $x^2 - 10x + 20$ \cr
 & $x^2 + 16x + 44$ \cr
 & $x^2 + 16x + 44$ \cr
\hline
\end{tabular}
\end{center}
\medskip
    \caption{Characteristic polynomials for Hecke operators on the
cuspidal subspace with norm level 3721\label{tab:3721-table}.  The
order corresponds to Table \ref{tab:heckereps}.}
\end{table}

\begin{table}[htb]
\begin{tabular}{|c|ccccc|}
\hline
$\Norm(\fn)$&$a_1$&$a_2$&$a_3$&$a_4$&$a_6$\cr
\hline
701 & $ -\z - 1$&$\z^2 - 1$&$1$&$-\z^2$&$0$\cr
2201 & $ -\z^2 - 2$&$\z^3 + \z^2$&$\z$&$0$&$0$\cr
2351 & $ 1$&$\z^2 + \z + 2$&$\z$&$\z^2 + 1$&$0$\cr
3061 & $ 2\z^3 + \z + 2$&$1$&$-\z^2$&$0$&$0$\cr
3355 & $ \z^3 - \z + 1$&$-\z$&$0$&$1$&$0$\cr
3571 & $ -\z^3 - \z$&$\z - 1$&$\z^3 + 1$&$0$&$0$\cr
3641a & $ -2\z^2 - \z - 1$&$\z^2$&$\z + 1$&$0$&$0$\cr
3641b & $ \z^3 - 1$&$2\z^3 + \z + 2$&$\z^2 + 1$&$\z$&$0$\cr
4681 & $ \z^3$&$-\z^3 + 1$&$\z^3$&$-\z^3$&$0$\cr
5081 & $ -\z^2 + \z + 1$&$\z^3 + \z + 1$&$\z + 1$&$0$&$0$\cr
5101 & $ -\z^3 - 2\z$&$-1$&$\z$&$0$&$0$\cr
6961 & $ \z^3 - 1$&$-\z - 2$&$0$&$\z + 1$&$0$\cr
\hline
\end{tabular}
\medskip
\caption{Equations for elliptic curves over $F$\label{tab:ec}.  The curve 3641b corresponds to the first eigenclass labelled 3641b in Table \ref{tab:qclasses}.  (We did not find a curve corresponding to the second eigenclass labelled 3641b.)}
\end{table}

\def\fv{\sqrt{5}}
\begin{table}
\begin{center}
\begin{tabular}{|p{15pt}|c|p{15pt}|c|p{15pt}|c|p{15pt}|c|}
\hline
$q $& $a_\fq$ & $q$ & $a_\fq$& $q$ & $a_\fq$& $q$ & $a_\fq$\cr
\hline
 2 &$\fv$&11 &$1$&29&$4\fv$&41&$-8$\cr
 3 &$-2\fv$&11 &$2$&29&$-2\fv $&41&$2$\cr
 5 &0&19 &$-2\fv $&31&$-8$&59&$-2\fv$\cr
 7 &$2\fv $&19 &$0$&31&$2$&59&$4\fv$\cr
\hline
\end{tabular}
\end{center}
\medskip
    \caption{Hecke eigenvalues of the Hilbert modular newform $g$
corresponding to the third eigenclass at $3025$.
\label{tab:3025-table}}
\end{table}

\begin{table}[h]
\begin{center}
\begin{tabular}
{|@{\hskip3pt}c@{\hskip3pt} |@{\hskip3pt}c@{\hskip3pt}
 |@{\hskip3pt}c@{\hskip3pt} |@{\hskip3pt}c@{\hskip3pt}|
 |@{\hskip3pt}c@{\hskip3pt} |@{\hskip3pt}c@{\hskip3pt}
 |@{\hskip3pt}c@{\hskip3pt} |@{\hskip3pt}c@{\hskip3pt}|
 |@{\hskip3pt}c@{\hskip3pt} |@{\hskip3pt}c@{\hskip3pt}
 |@{\hskip3pt}c@{\hskip3pt} |@{\hskip3pt}c@{\hskip3pt}|
 |@{\hskip3pt}c@{\hskip3pt} |@{\hskip3pt}c@{\hskip3pt}
 |@{\hskip3pt}c@{\hskip3pt} |@{\hskip3pt}c@{\hskip3pt}|}\hline
$(a,b,c,d)$&$s$&$g$&$I$&$(a,b,c,d)$&$s$&$g$&$I$&
$(a,b,c,d)$&$s$&$g$&$I$&$(a,b,c,d)$&$s$&$g$&$I$\\\hline
(0,0,1,1)&3&{\bf 3}&12&
(1,0,1,1)&3&1&2&
(0,0,2,1)&3&0&-&
(1,0,2,1)&1&0&-\\\hline
(0,1,1,1)&0&-&-&
(1,1,1,1)&2&0&-&
(0,1,2,1)&0&-&-&
(1,1,2,1)&2&1&2\\\hline
(0,2,1,1)&3&{\bf 3}&8&
(1,2,1,1)&1&0&-&
(0,2,2,1)&2&0&-&
(1,2,2,1)&2&1&2\\\hline
(0,3,1,1)&1&0&-&
(1,3,1,1)&1&0&-&
(0,3,2,1)&3&1&2&
(1,3,2,1)&1&{\bf 1}&{\bf 4}\\\hline
(0,4,1,1)&2&0&-&
(1,4,1,1)&0&-&-&
(0,4,2,1)&2&1&2&
(1,4,2,1)&2&1&2\\\hline
(0,5,1,1)&1&0&-&
(1,5,1,1)&1&0&-&
(0,5,2,1)&2&0&-&
(1,5,2,1)&4&{\bf 4}&{\bf 22}\\\hline
\end{tabular}
\end{center}
\medskip
\caption{Data concerning Mordell--Weil groups
of elliptic curves~$E(\Delta)$\label{table-app}}
\end{table}

\end{document}